\documentclass[11pt]{article}
\usepackage[a4paper]{anysize}\marginsize{3.5cm}{3.5cm}{1.3cm}{2cm}
\pdfpagewidth=\paperwidth \pdfpageheight=\paperheight
\usepackage{amsfonts,amssymb,amsthm,amsmath,eucal}
\usepackage{pgf}
\usepackage{bbm}

\usepackage{tikz} 
\usepackage{float}
\usetikzlibrary{arrows}
\usepackage{subfigure}
\usepackage{caption}

\theoremstyle{plain}
\newtheorem{thm}{Theorem}[section]
\newtheorem{theorem}[thm]{Theorem}

\newtheorem{proposition}[thm]{Proposition}

\theoremstyle{definition}
\newtheorem{definition}[thm]{Definition}
\newtheorem{remark}[thm]{Remark}
\newtheorem{example}[thm]{Example}

\newtheorem{thevarthm}[thm]{\varthmname}

\newenvironment{varthm*}[1]{\trivlist\item[]{\bf #1.}\it}{\endtrivlist}

\renewcommand\geq{\geqslant}

\renewcommand\leq{\leqslant}

\newcommand\be{\begin{eqnarray*}}
\newcommand\ee{\end{eqnarray*}}

\newcommand\Z{\mathbb Z}

\renewcommand\P{\mathbb P}

\newcommand\newop[2]{\def#1{\mathop{\rm #2}\nolimits}}
\newop\log{log}
\newop\ord{ord}
\newop\Gal{Gal}
\newop\SL{SL}
\newop\Bl{Bl}
\newop\mult{mult}
\newop\mass{mass}
\newop\div{div}
\newop\codim{codim}
\newop\sing{sing}
\newop\Zeroes{Zeroes}

\newcommand\wtilde[1]{\widetilde{#1}}

\newcommand\alphahat{\widehat{\alpha}}

\def\keywordname{{\bfseries Keywords}}%
\def\keywords#1{\par\addvspace\medskipamount{\rightskip=0pt plus1cm
\def\and{\ifhmode\unskip\nobreak\fi\ $\cdot$
}\noindent\keywordname\enspace\ignorespaces#1\par}}
\def\subclassname{{\bfseries Mathematics Subject Classification
(2000)}\enspace}
\def\subclass#1{\par\addvspace\medskipamount{\rightskip=0pt plus1cm
\def\and{\ifhmode\unskip\nobreak\fi\ $\cdot$
}\noindent\subclassname\ignorespaces#1\par}}

\newcommand\beginproof[1]{\trivlist\item[\hskip\labelsep{\em #1.}]}
\newcommand\proofof[1]{\beginproof{Proof of #1}}
\definecolor{uuuuuu}{rgb}{0,0,1}
\definecolor{qqqqff}{rgb}{0,0,1}
\definecolor{xdxdff}{rgb}{0,0,1}

\def\endproof{\hspace*{\fill}\endproofsymbol\endtrivlist}

\def\endproofsymbol{\frame{\rule[0pt]{0pt}{6pt}\rule[0pt]{6pt}{0pt}}}

\begin{document}

\author{M.~Dumnicki, T.~Szemberg, H.~Tutaj-Gasi\'nska}
\title{Symbolic powers of planar point configurations II}
\date{\today}
\maketitle
\thispagestyle{empty}

\begin{abstract}
   In \cite{DST13} we began to study the initial sequences $\alpha(I^{(m)})$, $m=1,2,3\dots,$
   of radical ideals $I$ of finite sets of points in the projective plane.
   In the present note we complete results obtained in \cite{DST13}
   by answering a number of questions left open in the previous note
   and we extend our considerations to the asymptotic setting of Waldschmidt constants.
   The concept of the Bezout decomposition introduced in Definition \ref{def:bezout decomposition}
   might be of independent interest.
\subclass{MSC 14C20 \and MSC 14J26 \and MSC 14N20 \and MSC 13A15 \and MSC 13F20}
\end{abstract}


\section{Introduction}
   Symbolic powers of ideals of point configurations have attracted
   considerable attention in recent years. Apart from traditional
   paths of research
   motivated by various problems in Algebraic Geometry, Commutative Algebra
   and Combinatorics (see e.g. \cite{BoHa10a}, \cite{HH07}, \cite{Sul08}),
   ideals of planar points have been recently studied in connection
   with the counterexamples to the $I^{(3)}\subset I^2$ containment
   (see e.g. \cite{DST13b}, \cite{BCH14}, \cite{DHNSST14}, \cite{HaSe13})
   and with the Bounded Negativity Conjecture (see \cite{BNAL}).
   These recent directions of investigation focus on special configurations
   of points. The study of such configurations, from a yet slightly
   different point of view, has been initiated by
   Bocci and Chiantini in \cite{BocCha11}. We follow their approach in the present note.

   For a homogeneous ideal $I=\bigoplus_{d\geq 0}I_d$, we define
   the \emph{initial degree} $\alpha(I)$ of $I$ as the least
   integer $d$ such that $I_d\neq 0$. More generally, we define
   the \emph{initial sequence} (or simply the $\alpha$-sequence of $I$)
   as the strictly increasing sequence of integers
   \begin{equation}\label{eq:initial sequence}
      \alpha(I)<\alpha(I^{(2)})<\alpha(I^{(3)})<\alpha(I^{(4)})<\ldots,
   \end{equation}
   where $I^{(m)}$ denotes the $m$-th symbolic power of $I$.

   There is a related asymptotic quantity introduced by Chudnovsky \cite{Chu81}
   and rediscovered by Harbourne who named it
   the \emph{Waldschmidt constant} of $I$:
   $$\alphahat(I)=\lim\limits_{m\to\infty}\frac{\alpha(I^{(m)})}{m}=\inf\limits_{m\geq 1}\frac{\alpha(I^{(m)})}{m}.$$
   For radical ideals $I=I(Z)$ attached to configurations $Z$ of points in $\P^2$
   Bocci and Chiantini
   studied the question to what extent the value of the difference $\alpha(I^{(2)})-\alpha(I)$
   determines the geometry of $Z$.
   Their results, still for planar points configurations, have been considerably
   generalized in \cite{DST13} and further extended to other varieties
   in \cite{BaczPhD}, \cite{P1P1}, \cite{BauSze15} and \cite{DLS15}.
   Many ideas presented here can be adapted to a more general setting.
   We don't dwell on this point in order to keep the presentation as
   transparent as possible.

   It is convenient to define the \emph{first differences sequence} of the initial
   sequence as $\beta_m(I)=\alpha(I^{(m+1)})-\alpha(I^{(m)})$ for $m\geq 1$ and to set
   $\beta_0(I)=\alpha(I)$. We call this
   sequence simply the $\beta$-sequence of $I$. Of course, the $\alpha$ sequence determines
   the $\beta$ sequence and vice versa.
   In the present note we focus on zero sets of ideals whose $\beta$-sequence contains
   many twos and threes.

   Our main result
   is the following classification statement,
   see Definition \ref{def:star config} for the terminology applied.
\begin{varthm*}{Main Theorem}
   Let $Z=\left\{P_1,\ldots,P_s\right\}$
   be a finite set of points in $\P^2$ and let $I=I(Z)$ be its radical ideal.
   If
   $$\alphahat(I)\;<\;\frac94,$$
   then $Z$
   \begin{itemize}
   \item[a)] is contained in a line ($\alphahat(I)=1$) or a conic ($\alphahat(I)\leq 2$) or;
   \item[b)] is a $4$-star ($\alphahat(Z)=2$) and $s=6$ in this case.
   \end{itemize}
   Moreover if $\alpha(I^{(m)})=9/4$ for some $m$, then $Z$ is a $3$-quasi star.
\end{varthm*}
   As sample consequences of the Main Theorem we derive the following results.
\begin{varthm*}{Corollary A}
   Let $I$ be the radical ideal of a finite set $Z$ of points in $\P^2$.
   \begin{itemize}
   \item[a)] If there exists $m\geq 1$ such that $\beta_m(I)=\beta_{m+1}(I)=1$,
      then $\alpha(I)=1$, i.e. the set $Z$ is contained in a line.
   \item[b)] If there exists $m\geq 1$ such that
      $\beta_m(I)=\beta_{m+1}(I)=\beta_{m+2}(I)=\beta_{m+3}(I)=\beta_{m+4}(I)=2$,
      then $\alpha(I)=2$ i.e. $Z$ is contained in a conic.
   \item[c)] For any $d\geq 3$ there exist configurations of points such
      that $\beta_m(I)=d$ for all $m\geq 2$ but $\alpha(I)\geq d+1$.
   \end{itemize}
\end{varthm*}
\begin{varthm*}{Corollary B}
   Let $I$ be the radical ideal of a finite set $Z$ of points in $\P^2$ with an integral Waldschmidt constant $\alphahat(I)=d$.
   \begin{itemize}
   \item[(i)] if $d=1$, then $\alpha(I)=1$, i.e. $Z$ is contained in a line;
   \item[(ii)] if $d\geq 2$, then $\alpha(I)$ need not to be equal $d$, i.e. there exist configurations of points $Z$ with $\alpha(I)\geq d+1$.
   \end{itemize}
\end{varthm*}
   We work over an algebraically closed field of characteristic $0$.

\section{Preliminaries and auxiliary results}
   In this section we fix the notation and recall very useful results
   of Chudnovsky and Esnault and Viehweg.

   For a point $P\in\P^2$ let $I(P)$ denote the radical ideal containing all forms
   vanishing at $P$.

   Let $Z=\left\{P_1,\dots,P_s\right\}$
   be a fixed finite set of points in the projective plane. Then the ideal of $Z$ is
   $$I=I(Z)=I(P_1)\cap\dots\cap I(P_s).$$
   In this situation, for a positive integer $m$, the $m$--th symbolic power of $I$ is defined by
   \begin{equation}\label{eq:symb power}
      I^{(m)}:=I(P_1)^m\cap\dots\cap I(P_s)^m,
   \end{equation}
   see \cite[Chapter IV.12, Definition]{ZarSam75} for the general definition.
   It is convenient and customarily to denote by $mZ$ the subscheme of $\P^2$ defined
   by the ideal $I^{(m)}$. We will also write $\alpha(mZ)$ rather than $\alpha(I^{(m)})$
   if we are primarily interested in the geometry of the set $Z$.

   The values of the initial sequence \eqref{eq:initial sequence} were considered
   already by Chudnovsky in \cite{Chu81} with the notation $\alpha(mZ)=\Omega(Z,m)$.
   He showed \cite[General Theorem 6]{Chu81} that one has always the following
   inequality (for sets $Z$ of points in $\P^2$):
   $$\frac{\alpha(Z)+1}{2}\leq \frac{\alpha(mZ)}{m}.$$
   This result has been generalized by Esnault and Viehweg \cite[Inequality (A), page 76]{EV83}.
   For any $n\leq m$ we have
   \begin{equation}\label{eq:EV}
   \frac{\alpha(nZ)+1}{n+1}\leq \frac{\alpha(mZ)}{m}.
   \end{equation}
   As the immediate corollary we have
   \begin{equation}\label{eq:EV2}
      \frac{\alpha(mZ)+1}{m+1}\leq\; \alphahat(Z)\leq\; \frac{\alpha(mZ)}{m}
   \end{equation}
   for any $m\geq 1$.
   We have also the following useful reformulation of \eqref{eq:EV}.
\begin{proposition}\label{prop:EV new}
   Let $I$ be the radical ideal of a finite set of points $Z$ in the projective plane.
   Let $d\geq k\geq 2$ and $m\geq 1$ be fixed integers such that
   \begin{equation}\label{eq:difference}
      \alpha((m+k)Z)=\alpha(mZ)+d.
   \end{equation}
   Then
   \begin{equation}\label{eq:EV new}
      \alpha((m+k)Z)\leq\frac{d-1}{k-1}(m+k).
   \end{equation}
\end{proposition}
\proof
   Plugging \eqref{eq:difference} into \eqref{eq:EV} we get
   $$\frac{\alpha((m+k)Z)-(d-1)}{m+1}\;\leq\;\frac{\alpha((m+k)Z)}{m+k}.$$
   Resolving this inequality with respect to $\alpha((m+k)Z)$ yields \eqref{eq:EV new}.
\endproof

   In the sequel we will encounter some interesting geometrical configurations of points.
\begin{definition}[Star configuration of points]\label{def:star config}
   We say that $Z\subset\P^2$ is a \emph{star configuration} of degree $d$
   (or a $d$-star for short)
   if $Z$ consists of \textbf{all} intersection points
   of $d$ \textbf{general} lines in $\P^2$.
\end{definition}
   Star configurations can be defined much more generally and they pop
   up frequently in situation similar to those studied here. We refer
   to \cite{GHM13} for a very nice introduction to this circle of ideas.

   We will need also the following modification of Definition \ref{def:star config}.
\begin{definition}[Quasi star configuration of points]\label{def:quasi star config}
   We say that $Z\subset\P^2$ is a \emph{quasi star configuration} of degree $d$
   (or a $d$-quasi star for short)
   if $Z$ consists of all intersection points
   of $d$ general lines in $\P^2$ and additionally there
   is exactly one more point from $Z$ on each of the lines, moreover these
   additional points are not collinear.
\end{definition}
   Note that if the extra points in the above definition were
   collinear $Z$ would be a $(d+1)$-star. Note also that a $d$-star
   contains exactly $d(d-1)/2$ points and a $d$-quasi star
   contains exactly $d(d+1)/2$ points.
   The figure below depicts a $3$-quasi star.
\unitlength.15mm
\begin{figure}[H]
\centering
\begin{picture}(220,220)(0,0)
\put(20,0){\line(0,1){220}}   
\put(0,20){\line(1,0){220}}   
\put(0,220){\line(1,-1){220}}   
 \put(20, 200){\circle*{7}} \put(30,200){C}
 \put(20, 80){\circle*{7}} \put(30,80){D}
 \put(20, 20){\circle*{7}} \put(30,30){A}
 \put(140, 20){\circle*{7}} \put(130,30){E}
 \put(200, 20){\circle*{7}} \put(210,30){B}
 \put(80, 140){\circle*{7}} \put(90,150){F}
\end{picture}
   \caption{}\label{fig:3-quasi star}
\end{figure}
   In the sequel we use without further comments
   the convention that the line passing through
   the points $P_i$ and $P_j$ is denoted by $L_{ij}$ and the line
   through $P$ and $Q$ by $L_{PQ}$.
\subsection{Bezout decomposition}
   We conclude this section with the following useful concept
   derived from the Bezout's Theorem, see \cite[Proposition 8.4]{Ful}.
   Let $Z=\left\{P_1,\ldots,P_s\right\}$ be a finite set of points in
   the projective plane $\P^2$. Let $D$ be an effective divisor of degree $d$
   vanishing to order $m_1,\ldots,m_s$ at the points $P_1,\ldots,P_s$
   respectively. Let $C_1,\ldots,C_r$ be irreducible curves of degrees
   $c_1,\ldots,c_r$ respectively and with $m^i_j=\mult_{P_j}C_i$
   for $i=1,\ldots,r$ and $j=1,\ldots,s$. Then we run the following
   algorithm. For each $i$ we compare the two numbers
   $$d_i:=D\cdot C_i=dc_i\;\mbox{ and }\; e_i:=\sum\limits_{j=1}^sm_jm^i_j.$$
   The first number is of course the intersection number of
   the divisor $D$ and the curve $C_i$. The second number appears
   in the Bezout's Theorem, which asserts that
   \begin{itemize}
      \item either $d_i\geq e_i$ holds;
      \item or $C_i$ is a component of $D$.
   \end{itemize}
   We decompose $D$ as
   $$D=\sum\limits_{i:\;e_i>d_i}C_i+D'$$
   and repeat the procedure for the divisor $D'$. After a finite number
   of steps we obtain the following decomposition
   \begin{equation}\label{eq:bezout}
      D=\sum_{i=1}^ra_iC_i+B(D)
   \end{equation}
   with $a_i\geq 0$ for $i=1,\ldots,r$ and the divisor $B(D)$ satisfying
   \begin{equation}\label{eq:bezout inequality}
      B(D)\cdot C_i\geq \sum_{j=1}^s\left(m_j-\sum_{k=1}^r a_km^k_j\right)m^i_j
   \end{equation}
   for all $i=1,\ldots,r$.
\begin{definition}[Bezout decomposition]\label{def:bezout decomposition}
   We call the decomposition in \eqref{eq:bezout} the \emph{Bezout decomposition}
   of $D$ with respect to the set $Z$ and curves $C_1,\ldots,C_r$
   and the divisor $B(D)$ the \emph{Bezout reduction} of $D$.
\end{definition}
\begin{remark}
   Note that the Bezout decomposition is determined purely numerically.
   That implies in particular that if $C_i$ and $C_{i'}$ are irreducible
   curves with the same degree and the same multiplicities in points from
   the set $Z$, then they will appear in \eqref{eq:bezout} with the
   same coefficients $a_i$ and $a_{i'}$.
\end{remark}
\begin{theorem}[Uniqueness of the Bezout decomposition]\label{thm:uniqueness of BD}
   The Bezout decomposition defined in \ref{def:bezout decomposition}
   is uniquely determined.
\end{theorem}
\proof
   Keeping the notation introduced in this paragraph let
   $f:X\to \P^2$ be the blow up of the points $P_1,\ldots,P_s$
   with exceptional divisors $E_1,\ldots,E_s$.
   We denote by $\wtilde{\Gamma}=f^*\Gamma-\sum_{i=1}^s\mult_{P_i}\Gamma\cdot E_i$
   the proper transform on $X$ of a curve $\Gamma\subset\P^2$.
   The condition $e_i>d_i$ is equivalent to $\wtilde{D}\cdot \wtilde{C_i}<0$.
   This immediately implies that $\wtilde{C_i}$ is a component of $\wtilde{D}$.
   Thus the Bezout decomposition of $D$ on $\P^2$ corresponds to subtracting from $\wtilde{D}$
   those curves among $\wtilde{C_1},\ldots,\wtilde{C_r}$ which have negative intersection
   with $\wtilde{D}$ (this can be viewed as a first reduction in taking the Zariski
   decomposition of $\wtilde{D}$). It suffices now to show that this reduction does not depend
   on the order in which the curves are subtracted.

   This is a consequence of the following simple observation. Suppose that $\wtilde{D}\cdot\wtilde{C}<0$
   and let $\wtilde{D}'$ be a divisor obtained from $\wtilde{D}$ by subtracting a curve $\wtilde{\Gamma}$
   different from $\wtilde{C}$. Then we have still $\wtilde{D}'\cdot\wtilde{C}<0$ and have to
   subtract the curve $\wtilde{C}$ from $\wtilde{D}'$ according to our algorithm.
   This shows that \emph{locally} the change of order in the reduction
   procedure does not influence the resulting divisor. This means in turn that the Bezout decomposition
   is \emph{locally} confluent. Since the algorithm has to stop after finitely many steps,
   it is also \emph{globally} confluent and the uniqueness of the Bezout reduction divisor $B(D)$
   follows easily from an elementary version of the Church-Rosser Theorem \cite{ChuRos36}.
\endproof
\section{Configurations of points with $\alphahat(Z)<\frac94$}
   In this section we will prove the Main Theorem.
   Of course if $Z$ is contained in a line or in a conic then
   $\alphahat(Z)=1$ or $\alphahat(Z)\leq2$ respectively (see Proposition \ref{prop:Waldschmidt less or equal 2}
   for the complete list of sets $Z$ with Waldschmidt constants $<2$). So we restrict our attention
   to sets not contained in a conic. Thus let $Z=\left\{P_1,\ldots,P_s\right\}$
   be a finite set of points in $\P^2$ not contained
   in a conic, so that in particular $s\geq 6$ holds. We can assume,
   renumbering the points if necessary, that the subset $W=\left\{P_1,\ldots,P_6\right\}$
   is not contained in a conic. In order to complete the proof
   of the Main Theorem we need to show that
   $Z=W$ and $Z$ is a $4$-star.
   The proof splits into several cases.

   Note to begin with that the assumption $\alphahat(Z)<\frac94$
   implies that there exists $m\geq 1$ such that
   \begin{equation}\label{eq:curve of deg leq 94m}
      \frac{\alpha(mZ)}{m}\leq\frac94  
   \end{equation}
   We fix such $m$ and write it in the form
   $$m=4n+p\;\mbox{ with }\; 0 \leq p \leq 3.$$
   Then
   $$\alpha(mZ)\leq 9n+2p.$$
   By \eqref{eq:curve of deg leq 94m} there exists a divisor $\Gamma$ of degree
   $9n+2p$ vanishing along $(4n+p)W$. We work with this divisor throughout the proof.
   We split the proof in a number of cases.

   {\bf Case 1.} Assume that no three points in $W$ are collinear. Then each conic, passing through
   exactly five points of $W$, is irreducible. We denote by $C_i$ the conic passing through all points
   in $W$ but $P_i$ for $i=1,\ldots,6$ and perform the Bezout decomposition
   of $\Gamma$ with respect to $W$ and the conics $C_1,\ldots,C_6$. Since the situation
   is symmetric with respect to these curves, we end up with the following decomposition
   $$\Gamma=k(C_1+\ldots+C_6)+B(\Gamma)$$
   with the inequality
   $$2(9n+2p-12k) \geq 5(4n+p-5k)$$
   coming from \eqref{eq:bezout inequality}. Equivalently we have
   \begin{equation}\label{eq:case 1}
      k \geq 2n+p.
   \end{equation}
   On the other hand the degree of the Bezout reduction divisor $B(\Gamma)$ must satisfy
   \begin{equation}\label{eq:case 1 2}
      9n+2p-12k \geq 0.
   \end{equation}
   But \eqref{eq:case 1} and \eqref{eq:case 1 2} are contradictory for $m>0$.
   Hence there are at least 3 collinear points in $W$.

   {\bf Case 2.} Three points in $W$ are collinear. Without loss of generality let $P_1, P_2, P_3$
   lie on a line $L$. We keep this assumption until the end of the proof.

   Observe that no more points of $W$ lie on $L$, since otherwise $W$ would lie on a
   conic. Similarly, $P_4, P_5, P_6$ cannot be collinear hence they determine three distinct lines
   $L_{45}, L_{46}$ and $L_{56}$.
   The union of these lines (the triangle determined by the set $\{P_4,P_5,P_6\}$) is denoted by $T=L_{45}+L_{46}+L_{56}$.
   $T$ may or may not contain some of the points points $P_1, P_2, P_3$. The rest of the proof
   splits onto subcases
   depending on how many of the points $P_1, P_2, P_3$ lie on $T$.

   {\bf Subcase 2.0.} The triangle does not pass through any of the points $P_1,P_2,P_3$.
   This situation is depicted below.
\begin{figure}[H]
\centering
   \begin{minipage}{0.4\textwidth}
   \centering
\begin{tikzpicture}[line cap=round,line join=round,>=triangle 45,x=1.0cm,y=1.0cm]
\clip(-2.4,0.0) rectangle (3.0,4.0); \draw [domain=-2.4:3.0] plot(\x,{(--10.3296--1.1*\x)/3.44});
\draw [domain=-2.4:3.0] plot(\x,{(--1.0736-1.42*\x)/1.8599999999999999}); \draw [domain=-2.4:3.0]
plot(\x,{(--0.37-1.16*\x)/-0.18}); \draw [domain=-2.4:3.0]
plot(\x,{(--1.7047999999999999--1.04*\x)/1.08});
\begin{scriptsize}
\draw [fill=black] (-1.76,2.44) circle (1.5pt); \draw[color=black] (-1.56,2.719999999999999) node
{$P_1$};
\draw [fill=black] (1.68,3.54) circle (1.5pt); \draw[color=black] (1.88,3.8199999999999994)
node {$P_3$};
\draw [fill=black] (-0.2657008801251186,2.9178282069367354) circle (1.5pt);
\draw[color=black] (-0.06,3.199999999999999) node {$P_2$};
\draw [fill=black] (-0.58,1.02) circle
(1.5pt); \draw[color=black] (-0.6,1.33) node {$P_4$};
\draw [fill=black]
(0.6629729729729731,2.2169369369369374) circle (1.5pt); \draw[color=black]
(0.9,2.1) node {$P_6$};
\draw [fill=black] (0.3652610641471905,0.2983490800596719)
circle (1.5pt); \draw[color=black] (0.6,0.5) node {$P_5$};
\end{scriptsize}
\end{tikzpicture}
  \caption{}\label{fig: 2.0}
   \end{minipage}
\end{figure}
   We consider now the Bezout decomposition of $\Gamma$
   with respect to the set $W$ and the lines $L$, $L_{45}, L_{46}$ and $L_{56}$.
   By symmetry the Bezout decomposition of $\Gamma$ has the shape
   $$\Gamma=kL+\ell T+B(\Gamma).$$
   It is easy to check that $k\geq 1$ in this case.
   Note that the degree of $B(\Gamma)$ is $9n+2p-k-3\ell$
   and this divisor vanishes along $(4n+p-k)X+(4n+p-2\ell)Y$,
   where $X=\left\{P_1,P_2,P_3\right\}$ and
   $Y=\left\{P_4,P_5,P_6\right\}$.
   We have the following three inequalities:
\begin{align}
   \label{case2.0.1} 9n+2p-k-3\ell & \geq 3(4n+p-k) \quad \textrm{(the intersection $B(\Gamma)\cdot L$)} \\
   \label{case2.0.2} 9n+2p-k-3\ell & \geq 2(4n+p-2\ell) \quad \textrm{(the intersection of $B(\Gamma)$ with a line in $T$)} \\
   \label{case2.0.3} 9n+2p-k-3\ell & \geq 0 \quad \textrm{(the nonnegativity of the degree of $B(\Gamma)$).}
\end{align}
   The inequality \eqref{case2.0.2} gives $n\geq k-\ell$ and together with \eqref{case2.0.1} we obtain
   $$2k-3\ell-p \geq 3n \geq 3k-3\ell,$$
   which implies $k+p \leq 0$, a contradiction.

{\bf Subcase 2.1} The triangle $T$ contains exactly one of the points $\{P_1, P_2, P_3\}$.
   Without loss of generality we may assume that it is the point $P_1$
   and that it lies on the line $L_{45}$.
   This situation is depicted in the figure below.
\begin{figure}[H]
\centering
   \begin{minipage}{0.4\textwidth}
   \centering
\begin{tikzpicture}[line cap=round,line join=round,>=triangle 45,x=1.0cm,y=1.0cm]
\clip(-2.2,0.0) rectangle (3.5,4.3); \draw [domain=-2.2:3.5]
plot(\x,{(--11.4056--1.2400000000000002*\x)/3.7800000000000002}); \draw [domain=-2.2:3.5]
plot(\x,{(--2.1367999999999996-1.42*\x)/1.8599999999999999}); \draw [domain=-2.2:3.5]
plot(\x,{(--1.0872000000000004-1.1600000000000001*\x)/-0.17999999999999994}); \draw
[domain=-2.2:3.5] plot(\x,{(--1.6976--0.56*\x)/0.96});
\begin{scriptsize}
\draw [fill=black] (-1.74,2.46) circle (1.5pt); \draw[color=black] (-1.6,2.725) node
{$P_1$}; \draw [fill=black] (2.02,3.68) circle (1.5pt); \draw[color=black] (2.22,3.959999999999999)
node {$P_3$}; \draw [fill=black] (-0.13125135852394765,2.9742984961455834) circle (1.5pt);
\draw[color=black] (0.06,3.259999999999999) node {$P_2$}; \draw [fill=black] (-0.46,1.5) circle
(1.5pt); \draw[color=black] (-0.26,1.7799999999999985) node {$P_4$}; \draw [fill=black]
(1.332227488151659,2.5454660347551346) circle (1.5pt); \draw[color=black] (1.54,2.819999999999999)
node {$P_6$}; \draw [fill=uuuuuu] (0.9973545499751368,0.3873959887286589) circle (1.5pt);
\draw[color=uuuuuu] (1.2000000000000002,0.6599999999999979) node {$P_5$};
\end{scriptsize}
\end{tikzpicture}
  \caption{}\label{fig: 2.1}
   \end{minipage}
\end{figure}
   Now, the situation is symmetric with respect to the lines $L$ (through $P_1, P_2$ and $P_3$)
   and $L_{45}$ (which passes also through the point $P_1$). It is also symmetric with
   respect to the lines $L_{46}, L_{56}$ and $L_{26}$, $L_{36}$ (which are not indicated in the
   above picture). Performing the Bezout decomposition of $\Gamma$ with respect to these lines
   now, we obtain
   $$\Gamma=k(L+L_{45})+\ell(L_{46}+L_{56}+L_{26}+L_{36})+B(\Gamma).$$
   The curve $B(\Gamma)$ has degree $9n+2p-2k-4\ell$ and it vanishes to order
   $4n+p-2k$ at $P_1$, $4n+p-4\ell$ at $P_6$ and $4n+p-k-\ell$ at all other points of $W$.
   Hence we have the following inequalities
\begin{align}
   \label{case2.1.1} 9n+2p-2k-4\ell & \geq 2(4n+p-k-\ell)+(4n+p-2k) \\
   \label{case2.1.2} 9n+2p-2k-4\ell & \geq (4n+p-k-\ell)+(4n+p-4\ell) \\
   \label{case2.1.3} 9n+2p-2k-4\ell & \geq 0.
\end{align}
   The second inequality gives $n \geq k-\ell$ whereas the first one implies $2(k-\ell)\geq 3n+p$. Hence
   $2n \geq 2(k-\ell) \geq 3n+p$, which is absurd.

{\bf Subcase 2.2} Two points out of $\{P_1,P_2,P_3\}$ lie on $T$.
   Without loss of generality we may assume that $P_1$ lies on the line $L_{45}$, and $P_2$ lies on $L_{46}$.
   This configuration is indicated in the figure below.
\begin{figure}[H]
\centering
   \begin{minipage}{0.4\textwidth}
   \centering
\begin{tikzpicture}[line cap=round,line join=round,>=triangle 45,x=1.0cm,y=1.0cm]
\clip(-2.2,0.0) rectangle (5.3,5);
\draw [domain=-2.2:5.3] plot(\x,{(--11.4056--1.2400000000000002*\x)/3.7800000000000002}); 
\draw [domain=-2.2:5] plot(\x,{(--2.1367999999999996-1.42*\x)/1.8599999999999999}); 
\draw [domain=-2.2:3.5] plot(\x,{(--1.0872000000000004-1.1600000000000001*\x)/-0.17999999999999994}); 
\draw [domain=-2.2:5.3] plot(\x,{(--1.6976--0.56*\x)/0.96});
\begin{scriptsize}
\draw [fill=black] (-1.74,2.46) circle (1.5pt); \draw[color=black] (-1.6,2.725) node
{$P_1$}; \draw [fill=black] (2.02,3.68) circle (1.5pt); \draw[color=black] (2.22,3.959999999999999)
node {$P_3$};
\draw [fill=black] (4.89,4.64) circle (1.5pt);
\draw[color=black] (5,4.3) node {$P_2$};
\draw [fill=black] (-0.46,1.5) circle
(1.5pt); \draw[color=black] (-0.26,1.7799999999999985) node {$P_4$}; \draw [fill=black]
(1.332227488151659,2.5454660347551346) circle (1.5pt); \draw[color=black] (1.54,2.819999999999999)
node {$P_6$}; \draw [fill=uuuuuu] (0.9973545499751368,0.3873959887286589) circle (1.5pt);
\draw[color=uuuuuu] (1.2000000000000002,0.6599999999999979) node {$P_5$};
\end{scriptsize}
\end{tikzpicture}
  \caption{}\label{fig: 2.2}
   \end{minipage}
\end{figure}
   Let now $X=\{P_1,P_2,P_4\}$ and $Y=\{P_3,P_5,P_6\}$.
   Note that the set $Y$ lies on the triangle $T_X$ defined by $X$,
   each point on exactly one line. However, no point from $X$ lies on the triangle $T_Y$ defined by the set $Y$.
   If $Z = W$, then $Z$ is a $3$-quasi star. We show in Proposition \ref{prop:3-quasi star}
   that $\alphahat(Z)=9/4$ in this case.

   Thus we may assume that there exists an extra point $P_7\in Z$.
   Applying the Bezout decomposition with respect to the lines in triangles $T_X$ and $T_Y$
   we get
   $$\Gamma=kT_X+\ell T_Y+B(\Gamma)$$
   with $B(\Gamma)$ vanishing to order $4n+p-2k$ along $X$ and to order $4n+p-k-2\ell$
   along $Y$. Thus we obtain the following inequalities:
\begin{align}
   \label{case2.2.1} 9n+2p-3k-3\ell & \geq 2(4n+p-2k)+(4n+p-k-2\ell) \\
   \label{case2.2.2} 9n+2p-3k-3\ell & \geq 2(4n+p-k-2\ell)
\end{align}
   We observe additionally that removing each triangle $T_X$ or $T_Y$ from $\Gamma$
   causes the multiplicity of the residual divisor at the point $P_7$ to drop at most by one.
   So comparing the degree of the divisor $B(\Gamma)$ and its multiplicity at $P_7$
   we obtain the following inequality
\begin{align}
   \label{case2.2.3} 9n+2p-3k-3\ell & \geq 4n+p-k-\ell
\end{align}
   From \eqref{case2.2.1} we get $2k-\ell \geq 3m+p$, and from \eqref{case2.2.2} we get $\ell-k \geq -n$.
   Beginning with the inequality \eqref{case2.2.3} (after the reduction of terms) we obtain
   $$5n \geq 2k+2\ell-p = 6(\ell-k)+4(2k-\ell)-p \geq -6n+4(3n+p)-p=6n+3p$$
   which is absurd.

{\bf Subcase 2.3} The set $\{P_1,P_2,P_3\}$ lies on a triangle defined by $\{P_4,P_5,P_6\}$.
   Then $W$ is a $4$-star. If $Z=W$ we are done. Otherwise consider an extra point $P_7\in Z$.
   We denote by $\Delta$ the union of the $4$ lines determined by $W$.
   Our two final cases depend on whether $P_7$ lies on $\Delta$ or not.

{\bf Subsubcase 2.3.a} We assume first that the point $P_7$ does not lie on $\Delta$.
   The situation is indicated in the figure below.
\begin{figure}[H]
\centering
   \begin{minipage}{0.4\textwidth}
   \centering
\begin{tikzpicture}[line cap=round,line join=round,>=triangle 45,x=1.0cm,y=1.0cm]
\clip(-2.2,1.0) rectangle (3.0,4.0);
\draw [domain=-2.2:3.0] plot(\x,{(--10.3296--1.1*\x)/3.44});
\draw [domain=-2.2:3.0] plot(\x,{(--2.0252-1.42*\x)/1.8599999999999999});
\draw [domain=-2.2:3.0] plot(\x,{(-1.0987349829394668-1.8339744723872107*\x)/-0.20224744092000524});
\draw [domain=-2.2:3.0] plot(\x,{(-7.7004-1.78*\x)/-3.02});
\begin{scriptsize}
\draw [fill=black] (-1.76,2.44) circle (1.5pt); \draw[color=black] (-1.56,2.719999999999999) node
{$P_1$}; \draw [fill=black] (1.68,3.54) circle (1.5pt); \draw[color=black]
(1.88,3.8199999999999994) node {$P_3$}; \draw [fill=black]
(-0.27775255907999474,2.9139744723872107) circle (1.5pt); \draw[color=black]
(-0.08000000000000002,3.199999999999999) node {$P_2$}; \draw [fill=black] (-0.46,1.44) circle
(1.5pt); \draw[color=black] (-0.26,1.7199999999999984) node {$P_4$}; \draw [fill=black]
(-0.3400132501799255,2.3493961638012366) circle (1.5pt); \draw[color=black]
(-0.14,2.6199999999999988) node {$P_6$}; \draw [fill=uuuuuu] (-1.079934729971576,1.913283503526687)
circle (1.5pt); \draw[color=uuuuuu] (-0.88,2.1999999999999984) node {$P_5$}; \draw [fill=black]
(1.46,2.32) circle (1.5pt); \draw[color=black] (1.6600000000000001,2.5999999999999988) node
{$P_7$};
\end{scriptsize}
\end{tikzpicture}
   \caption{}\label{fig: 2.3.a}
   \end{minipage}
\end{figure}
   Let
   $$\Gamma=k\Delta+B(\Gamma)$$
   be the Bezout decomposition of $\Gamma$ with respect to the lines in $\Delta$,
   so that
\begin{align}
   \label{case2.3.a.1} 9n+2p-4k & \geq 3(4n+p-2k)
\end{align}
   holds. Comparing the degree $9m+2p-4k$ of the divisor $B(\Gamma)$
   with its multiplicity $4n+p$ at $P_7$ we obtain additionally that
\begin{align}
   \label{case2.3.a.2} 9n+2p-4k & \geq 4n+p
\end{align}
   Reducing terms in inequalities \eqref{case2.3.a.1} and \eqref{case2.3.a.2} we get
   $$2k\geq 3n+p\;\;\mbox{ and }\;\; 5n\geq 4k-p,$$
   which gives a contradiction.

{\bf Subsubcase 2.3.b}
   Now we pass to the final case and assume that $P_7$ is contained in $\Delta$.
   Without loss of generality we assume that $P_7$ lies on the line $L$
   defined by $P_1$ and $P_2$. This is indicated by the following figure.
\begin{figure}[H]
\centering
   \begin{minipage}{0.4\textwidth}
   \centering
\begin{tikzpicture}[line cap=round,line join=round,>=triangle 45,x=1.0cm,y=1.0cm]
\clip(-2.2,1.0) rectangle (3.6,4.7);
\draw [domain=-2.2:3.6] plot(\x,{(--10.3296--1.1*\x)/3.44});
\draw [domain=-2.2:3.6] plot(\x,{(--2.0252-1.42*\x)/1.8599999999999999});
\draw [domain=-2.2:3.6] plot(\x,{(-1.0987349829394668-1.8339744723872107*\x)/-0.20224744092000524});
\draw [domain=-2.2:3.6] plot(\x,{(-7.7004-1.78*\x)/-3.02});
\begin{scriptsize}
\draw [fill=black] (-1.76,2.44) circle (1.5pt); \draw[color=black] (-1.56,2.7199999999999993) node
{$P_1$}; \draw [fill=black] (1.68,3.54) circle (1.5pt); \draw[color=black] (1.88,3.82) node
{$P_3$}; \draw [fill=black] (-0.27775255907999474,2.9139744723872107) circle (1.5pt);
\draw[color=black] (-0.08000000000000002,3.1999999999999997) node {$P_2$}; \draw [fill=black]
(-0.46,1.44) circle (1.5pt); \draw[color=black] (-0.26,1.7199999999999993) node {$P_4$}; \draw
[fill=black] (-0.3400132501799255,2.3493961638012366) circle (1.5pt); \draw[color=black]
(-0.14,2.6199999999999997) node {$P_6$}; \draw [fill=uuuuuu] (-1.079934729971576,1.913283503526687)
circle (1.5pt); \draw[color=uuuuuu] (-0.88,2.1999999999999993) node {$P_5$}; \draw [fill=black]
(2.76,3.86) circle (1.5pt); \draw[color=black] (2.96,4.14) node {$P_7$};
\end{scriptsize}
\end{tikzpicture}
   \caption{}\label{fig: 2.3.b}
   \end{minipage}
\end{figure}
    Let now $X=\{P_1,P_2,P_3\}$ and $Y=\{P_4,P_5,P_6\}$.
    For the Bezout decomposition we consider now the following divisors:
    $L$, the triangle $T_Y$ and the pencil $\Pi=L_{74}+L_{75}+L_{76}$.
    Let
    $$\Gamma=kL+\ell T_Y+t\Pi+B(\Gamma)$$
    be the Bezout decomposition.
    The divisor $B(\Gamma)$ has degree $9n+2p-k-3\ell-3t$ and vanishes
    along the points in $X$ to order $4n+p-k-\ell$, along $Y$ to order $4n+p-2\ell-t$
    and to order $4n+p-k-3t$ at $P_7$. So we have the following system of inequalities:
\begin{align}
   \label{case2.3.b.1} 9n+2p-k-3\ell-3t & \geq 3(4n+p-k-\ell)+(4n+p-k-3t)\\
   \label{case2.3.b.2} 9n+2p-k-3\ell-3t & \geq (4n+p-k-\ell)+2(4n+p-2\ell-t)\\
   \label{case2.3.b.3} 9n+2p-k-3\ell-3t & \geq (4n+p-2\ell-t)+(4n+p-k-3t)\\
   \label{case2.3.b.4} 9n+2p-k-3\ell-3t & \geq 0
\end{align}
   After reductions, the first three inequalities \eqref{case2.3.b.1}--\eqref{case2.3.b.3} give the following simpler system of inequalities:
\begin{align*}
   3k & \geq 7n+p, \\
   2\ell-t & \geq 3n+p, \\
   t-\ell & \geq -n.
\end{align*}
   Hence we have
   $$k\geq\frac13(7n+2p),\;\; \ell\geq 2n+p\;\mbox{ and }\; t\geq n+p,$$
   which combined with \eqref{case2.3.b.4} gives
   $$0\leq 9n+2p-k-3\ell-3t\leq 9n+2p-\frac13(7n+2p)-3(2n+p)-3(n+p) \leq -\frac73n-\frac{14}{3}p<0.$$
   Thus we are done with claims a) and b) of the Main Theorem.

   The ''moreover'' part follows from Case 2.2 above and the contradictions in all other cases.
\endproof
   Now we will compute the Waldschmidt constant of the $3$-quasi star.
\begin{proposition}[A $3$--quasi star]\label{prop:3-quasi star}
   Let $Z$ be a $3$--quasi star. Then $\alphahat(Z)=\frac94$.
\end{proposition}
\proof
   We use the notation as in the Figure \ref{fig:3-quasi star},
   thus let $Z=\left\{A,B,C,D,E,F\right\}$.
   Note to begin with that taking
   $$\Delta=2(L_{AB}+L_{BC}+L_{AC})+L_{DE}+L_{EF}+L_{DF}$$
   we obtain a divisor of degree $9$ vanishing to order $4$ along $Z$.
   This shows that
   $$\alphahat(Z)\leq\frac94.$$
   In order to prove the reverse inequality, assume that there exists
   a divisor $\Gamma$ of degree $d$ vanishing along $Z$ to order $m$
   and such that
   \begin{equation}\label{eq:less than 9/4}
      \frac{d}{m}<\frac94.
   \end{equation}
   We may also assume
   that $\Gamma$ has the least degree $d$ such that \eqref{eq:less than 9/4}
   holds.
   It is easy to check, using the Bezout's Theorem, that the divisor $\Delta$
   defined above has to be contained in $\Gamma$. But then, for
   $\Gamma'=\Gamma-\Delta$ we have $d'=d-9$ and $m'=m-4$.
   Hence $d'/m'<9/4$ holds as well, and this contradicts the minimality of $d$.
\endproof
\begin{remark}\label{rmk:9/4}
   We don't know if the above Proposition can be reversed, i.e. if a $3$-quasi star
   is the only configuration with the Waldschmidt constant equal $9/4$. There
   could exist a set $Z\subset\P^2$ such that $\alpha(mZ)>9/4$ for all $m\geq 1$
   but $\alphahat(Z)=9/4$.
\end{remark}
   We conclude this section with the classification of all point configurations
   with Waldschmidt constants less than $2$.
\begin{proposition}[Waldschmidt constants $<2$]\label{prop:Waldschmidt less or equal 2}
   Let $Z$ be a finite set of points in $\P^2$ with $\alphahat(Z)<2$.
   Then
   $$\alphahat(Z)=\frac{2k-1}{k}\;\;\mbox{ for some }\; k\in\Z_{>0}$$
   and $Z$ consists of $k$ points $P_1,\ldots,P_k$ contained in a line $L$
   and a single point $Q$ not contained in $L$.
\end{proposition}
\proof
   Suppose that $Z$ is not of the form asserted in the Proposition.
   Then there exist points $P,Q,R,S\in Z$ such that no $3$ of them are collinear.
   There exists a divisor $\Gamma$ of degree $d$ vanishing to order $m$
   along $Z$, so in particular at points $P,Q,R,S$, with $d<2m$.
   Let $M$ be the line through $P$ and $Q$ and let $N$ be the line determined by $R$ and $S$.
   Let
   $$\Gamma=k(M+N)+B(\Gamma)$$
   be the Bezout decomposition of $\Gamma$ with respect to the lines $M$ and $N$.
   Then $\deg(B(\Gamma))=d-2k$ and
   $$d-2k\geq 2(m-k)$$
   holds. But this is equivalent to $d\geq 2m> d$. A contradiction.

   Now we calculate $\alphahat(Z)$ for sets $Z$ described
   in the Proposition. Let $L_i$ be the line through $Q$ and $P_i$.
   For the divisor
   $$\Gamma=(k-1)L+L_1+\ldots+L_k$$
   we have $\deg(\Gamma)=2k-1$ and $\mult_P(\Gamma)=k$ for all $P\in Z$.
   Hence $\alphahat(Z)\leq \frac{2k-1}{k}$.

   It remains to show the Waldschmidt constant is exactly $\frac{2k-1}{k}$.
   To this end assume that there exists a divisor $\Gamma$ vanishing
   to order $m$ at all points in $Z$ and of degree $d$ satisfying
   $d<\frac{2k-1}{k}m$. Let $\Pi=L_1+\ldots+L_k$. Let
   $$\Gamma=pL+q\Pi+B(\Gamma)$$
   be the Bezout decomposition of $\Gamma$ with respect to the lines $L$ and $L_1,\ldots,L_k$.
   Then we have the following inequalities:
\begin{align*}
   d-p-kq & \geq k(m-p-q)\\
   d-p-kq & \geq m-p-q+m-kq.
\end{align*}
   Of course $\deg(B(\Gamma))=d-p-kq\geq 0$ holds. But all these three inequalities
   are contradictory.
\endproof
\section{Configurations with low $\beta$-sequences}
   In this section we prove Corollaries announced in the Introduction.
\proofof{Corollary A}
   Part a) was already proved as Corollary 3.5 in \cite{DST13}.

   Part b) with $m=1$ was also already proved as Theorem 4.11 in \cite{DST13}
   and the general case was conjectured there \cite[Conjecture 4.6]{DST13}.
   The argument presented here is new and covers both cases.

   We have by assumption that $\alpha((m+5)Z)=\alpha(mZ)+10$, hence Proposition \ref{prop:EV new}
   with $k=5$ and $d=10$ gives
   $$\alpha((m+5)Z)\leq\frac94(m+5).$$
   Thus the Main Theorem and Proposition \ref{prop:3-quasi star} applies.
   Since the $\beta$-sequence for the $4$-star is
   $$3,1,3,1,3,1,3,1,\ldots$$
   and for the $3$-quasi star it is
   $$3,2,2,2,3,2,2,2,3,2,2,2,\ldots$$
   these cases are excluded and we conclude that $Z$ is contained in a conic.

   Part c) for $d\geq 4$ follows from \cite[Example 4.14]{DST13}
   and for $d=3$ we construct a new example in the proof of Proposition \ref{prop:exp:4 2 3 3} below.
\endproof
   Now we pass to the second corollary.
\proofof{Corollary B}
   Part a) follows immediately from Chudnovsky inequality \eqref{eq:EV2} with $m=1$.

   Part b) for $d\geq 3$ follows from \cite[Example 4.14]{DST13},
   whereas for $d=2$ a $4$-star provides an example.
\endproof
\begin{proposition}\label{prop:exp:4 2 3 3}
   There exists a configuration $Z$ of points in $\P^2$ such that for the radical ideal $I=I(Z)$
   $$\alpha(I^{(m)})=3m\;\mbox{ for all }\; m\geq 2$$
   but $\alpha(I)=4$.
\end{proposition}
\proof
    We provide an explicit example. To this end we consider the subscheme $Z$
    consisting of ten points $P_1,\dots,P_{10}$ such that $P_1$, $P_2$, $P_3$ are the
    intersection points of three general lines $L_1,L_2,L_3$, points $P_4,\dots,P_9$ lie in pairs
    in general position
    on lines $L_1, L_2, L_3$ and
    $P_{10}$ is a general point in $\P^2$. This is indicated in the figure below.
\begin{figure}[H]
\centering
\begin{tikzpicture}[line cap=round,line join=round,>=triangle 45,x=1.0cm,y=1.0cm,scale=0.8]
\clip(0,1) rectangle (6.8,6.3);
\draw [domain=-4.3:15.8] plot(1.55,\x);
\draw [domain=-4.3:15.8] plot(\x,{(--5.344-0.0*\x)/3.3400000000000003});
\draw [domain=-4.3:15.8] plot(\x,{(--24.192-3.8400000000000003*\x)/3.3600000000000003});
\begin{scriptsize}
\draw [fill=qqqqff] (1.54,5.44) circle (1.5pt);
\draw[color=qqqqff] (1.7200000000000002,5.78) node {$P_3$};
\draw [fill=qqqqff] (1.56,1.6) circle (1.5pt);
\draw[color=qqqqff] (1.7400000000000002,1.94) node {$P_1$};
\draw [fill=qqqqff] (4.9,1.6) circle (1.5pt);
\draw[color=qqqqff] (5.08,1.94) node {$P_2$};
\draw [fill=xdxdff] (1.5468737284687373,4.120244134002442) circle (1.5pt);
\draw[color=xdxdff] (1.7200000000000002,4.46) node {$P_4$};
\draw [fill=xdxdff] (1.554894886748949,2.5801817442018176) circle (1.5pt);
\draw[color=xdxdff] (1.7400000000000002,2.9200000000000004) node {$P_5$};
\draw [fill=xdxdff] (2.58,1.6) circle (1.5pt);
\draw[color=xdxdff] (2.7600000000000002,1.94) node {$P_6$};
\draw [fill=xdxdff] (3.84,1.6) circle (1.5pt);
\draw[color=xdxdff] (4.02,1.94) node {$P_7$};
\draw [fill=xdxdff] (3.8109734513274343,2.8446017699115043) circle (1.5pt);
\draw[color=xdxdff] (4.0,3.18) node {$P_8$};
\draw [fill=xdxdff] (2.5484955752212395,4.287433628318584) circle (1.5pt);
\draw[color=xdxdff] (2.72,4.62) node {$P_9$};
\draw [fill=qqqqff] (5.16,4.76) circle (1.5pt);
\draw[color=qqqqff] (5.38,5.1) node {$P_{10}$};
\end{scriptsize}
\end{tikzpicture}
   \caption{}\label{fig:4 1}
\end{figure}
    It is easy to see that $\alpha(Z)=4$ in this situation.
    We claim that for $k\geq 2, \alpha(kZ)=3k$.  Indeed, $\alpha(kZ)\leq 3k$
    as for every $k\geq 2$
    there exists a divisor of degree $3k$ passing through $Z$ with multiplicities $k$ composed of:
\begin{enumerate}
   \item $t(L_1+L_2+L_3)+tC$, if $k=2t$
   \item $t(L_1+L_2+L_3)+(t-1)C+S$, if $k=2t+1$
\end{enumerate}
    where $C$ is a cubic with a double point in $P_{10}$ passing through $P_4,\dots,P_9$, while $S$ is
    the sextic passing through $P_1$, $P_2$, $P_3$, having multiplicity $2$ at $P_4,\dots,P_9$
    and multiplicity $3$ at $P_{10}$.

    It suffices to prove that there is no divisor of degree $3k-1$ passing through $Z$ with
    multiplicities $k$. Suppose that such a divisor $D$ exists.
    Let
    $D=a(L_1+L_2+L_3)+bC +B(D)$
    be the Bezout decomposition of $D$ with respect to the curves $L_1,L_2,L_3$ and $C$.
    As the degree of $B(D)$ is
    $3k-3a-3b-1$ and $B(D)$ passes through $Z$ with multiplicities\\
    \begin{center}
    $k-2a$ at $P_1,P_2,P_3$; $k-a-b$ at $P_4,\dots,P_9$ and $k-2b$ at $P_{10}$,
    \end{center}
    we obtain
    $$3k-3a-3b-1\geq 2(k-2a)+2(k-a-b)$$
    $$3(3k-3a-3b-1)\geq 2(k-2b)+6(k-a-b).$$
    Adding these inequalities and reducing terms we get $-4 \geq 0$, a contradiction.
\endproof
\begin{remark}\label{rmk:3 3 3 3}
   In the example above the $\beta$-sequence is $2,3,3,3,\ldots$ with $\beta_0=4$.
   We don't know if there exists a set $Z\subset\P^2$ with the $\beta$-sequence
   constantly equal $3$ and $\beta_0=4$. It would be desirable to know this
   in the view of \cite[Example 4.14]{DST13}.
\end{remark}
   Motivated by the above Remark, we provide now an example of a subscheme $Z$
   such that there is a considerable number of threes in the beginning of its
   $\beta$-sequence.
\begin{example}
   Let $\{P_1,\dots,P_7\}$ be generic points on an irreducible conic.
   Let $P_{8}$ be the intersection point
   of $L_{1,2}$ and $L_{6,7}$, let $P_{10}$ be the intersection point of $L_{2,3}$ and $L_{4,5}$
   and let $P_9$ be the intersection point of $L_{1,10}$ and
   $L_{6,7}$. This assumptions are illustrated in Figure \ref{fig:3 3 3 3}.
\definecolor{uuuuuu}{rgb}{0.0,0.0,1.0}
\definecolor{xdxdff}{rgb}{0.0,0.0,1.0}
\definecolor{qqqqff}{rgb}{0.0,0.0,1.0}
\begin{figure}[H]
\centering
   \begin{minipage}{0.4\textwidth}
   \centering
\begin{tikzpicture}[line cap=round,line join=round,>=triangle 45,x=1.0cm,y=1.0cm]
\clip(-2.3,-1.4) rectangle (5.9,3.9);
\draw [rotate around={-1.4975877390175145:(2.0499999999999994,0.84)}] (2.0499999999999994,0.84) ellipse (2.126386218328932cm and 1.4761498397856558cm);
\draw [domain=-4.3:15.799999999999997] plot(\x,{(--0.5579699575992105--0.1232487893693478*\x)/3.8404866783666662});
\draw [domain=-4.3:15.799999999999997] plot(\x,{(--0.48022584205858154--1.0117051216213613*\x)/0.5031153586326651});
\draw [domain=-4.3:15.799999999999997] plot(\x,{(--1.9580759244542936--0.47424456859620934*\x)/1.1894439641178578});
\draw [domain=-4.3:15.799999999999997] plot(\x,{(--5.568822412380808-0.9875577993774278*\x)/1.2564355737040511});
\draw [domain=-4.3:15.799999999999997] plot(\x,{(-2.0857954440833857-1.7782667720140812*\x)/-2.4257143477634644});
\begin{scriptsize}
\draw [fill=qqqqff] (2.84,2.2) circle (1.5pt);
\draw[color=qqqqff] (3.039999999999999,2.48) node {$P_4$};
\draw [fill=xdxdff] (1.618497214438312,2.29152403973163) circle (1.5pt);
\draw[color=xdxdff] (1.699999999999999,2.72) node {$P_3$};
\draw [fill=xdxdff] (0.42905325032045427,1.8172794711354205) circle (1.5pt);
\draw[color=xdxdff] (0.2799999999999994,2.16) node {$P_2$};
\draw [fill=xdxdff] (4.096435573704051,1.2124422006225724) circle (1.5pt);
\draw[color=xdxdff] (4.299999999999999,1.5) node {$P_5$};
\draw [fill=xdxdff] (4.029496801624967,0.274600759874439) circle (1.5pt);
\draw[color=xdxdff] (4.299999999999999,0.62) node {$P_6$};
\draw [fill=xdxdff] (-0.07406210831221083,0.8055743495140593) circle (1.5pt);
\draw[color=xdxdff] (0.3199999999999994,0.92) node {$P_1$};
\draw [fill=xdxdff] (0.1890101232583008,0.1513519705050912) circle (1.5pt);
\draw[color=xdxdff] (0.3799999999999994,0.44000000000000006) node {$P_7$};
\draw [fill=uuuuuu] (-0.40894615708638565,0.1321623745449188) circle (1.5pt);
\draw[color=uuuuuu] (-0.2000000000000005,-0.06) node {$P_8$};
\draw [fill=uuuuuu] (2.3516522394512536,2.5838411215281405) circle (1.5pt);
\draw[color=uuuuuu] (2.459999999999999,3.0) node {$P_{10}$};
\draw [fill=uuuuuu] (-1.0193785297132556,0.11257239618568715) circle (1.5pt);
\draw[color=uuuuuu] (-1.1400000000000006,0.4) node {$P_9$};
\end{scriptsize}
\end{tikzpicture}
   \caption{}\label{fig:3 3 3 3}
   \end{minipage}
\end{figure}
   Then using any algebra computer program (the authors did it with Singular \cite{DGPS})
   one can check that this configuration of  $\{P_1,\dots,P_{10}\}$ indeed
   satisfies $\alpha(Z)=4$,  $\alpha(kZ)-\alpha((k-1)Z)=3$ for $k=2,\dots,29$ and
   $\alpha(30Z)-\alpha(29Z)=4$.

\end{example}
   We conclude this note with the following.
\begin{remark}
   It would be interesting to see to what extend the classification
   stated in the Main Theorem can be prolonged for higher values of $\alphahat(Z)$.
   We expect the problem to be feasible for all $Z$ with $\alphahat(Z)<3$.
   We hope to come back to this issue soon.
\end{remark}
\paragraph*{Acknowledgement.} The final version of this paper was written during our
   stay at the Conference Center of the Polish Academy of Sciences in B\c edlewo.
   We would like to thank the Center for providing excellent working conditions
   and the Jagiellonian University in Cracow for providing financial support.
   We would like to thank Justyna Szpond for helpful remarks.


\bigskip \small

\bigskip
   Marcin Dumnicki,
   Jagiellonian University, Institute of Mathematics, {\L}ojasiewicza 6, PL-30348 Krak\'ow, Poland

\nopagebreak
   \textit{E-mail address:} \texttt{Marcin.Dumnicki@im.uj.edu.pl}

\bigskip
   Tomasz Szemberg,
   Department of Mathematics, Pedagogical University of Cracow,
   Podchor\c a\.zych 2,
   PL-30-084 Krak\'ow, Poland.

\nopagebreak
   \textit{E-mail address:} \texttt{tomasz.szemberg@gmail.com}

\bigskip
   Halszka Tutaj-Gasi\'nska,
   Jagiellonian University, Institute of Mathematics, {\L}ojasiewicza 6, PL-30348 Krak\'ow, Poland

\nopagebreak
   \textit{E-mail address:} \texttt{Halszka.Tutaj@im.uj.edu.pl}


\end{document}